\documentclass[11pt]{amsart}
\usepackage{amsmath,amssymb}
\newtheorem{theorem}{Theorem}[section]

\newtheorem{corollary}[theorem]{Corollary}

\theoremstyle{definition}

\begin{document}

\title{Minimal hypersurfaces and geometric inequalities}
\author{Simon Brendle}
\address{Department of Mathematics \\ Columbia University \\ New York NY 10027}
\begin{abstract}
In this expository paper, we discuss some of the main geometric inequalities for minimal hypersurfaces. These include the monotonicity formula, the Alexander-Osserman conjecture, the isoperimetric inequality for minimal surfaces, and the Michael-Simon Sobolev inequality. 
\end{abstract}
\thanks{The author was supported by the National Science Foundation under grant DMS-1806190 and by the Simons Foundation.}

\maketitle 

\section{Minimal surfaces} 

The study of minimal surfaces has a long history in geometry. By definition, a hypersurface in Euclidean space is minimal if it is critical point of the area functional. Alternatively, minimal surfaces can be characterized in terms of their extrinsic curvature. To explain this, suppose $\Sigma$ is a hypersurface in $\mathbb{R}^{n+1}$ (possibly with boundary $\partial \Sigma$). Moreover, let $F_s: \Sigma \to \mathbb{R}^{n+1}$ denote a one-parameter family of immersions with $F_0(x)=x$, and let $V(x) := \frac{\partial}{\partial s} F_s(x) \Big |_{s=0}$ denote the velocity of the deformation. Let $\Sigma_s := F_s(\Sigma)$. If $V$ vanishes along the boundary $\partial \Sigma$, then the first variation of area is given by 
\[\frac{d}{ds} |\Sigma_s| \Big |_{s=0} = \int_\Sigma H \, \langle V,\nu \rangle.\] 
Here, $\nu$ denotes the unit normal to $\Sigma$ and $H$ denotes the scalar mean curvature. In other words, the mean curvature vector is $-H \, \nu$ in our convention.

We can think of the mean curvature as follows. Given two tangential vector fields $X$ and $Y$ along the hypersurface $\Sigma$, the second fundamental form is defined by $h(X,Y) = \langle \bar{D}_X \nu,Y \rangle$, where $\bar{D}$ denotes the standard flat connection on $\mathbb{R}^{n+1}$. The second fundamental form is a symmetric bilinear form on the tangent space to $\Sigma$, i.e. $h(Y,X) = h(X,Y)$. The mean curvature is defined as the trace of the second fundamental form; that is, 
\[H = \sum_{i=1}^n h(e_i,e_i) = \sum_{i=1}^n \langle \bar{D}_{e_i} \nu,e_i \rangle,\] 
where $\{e_1,\hdots,e_n\}$ is a local orthonormal frame on the hypersurface $\Sigma$. In view of the discussion above, we can think of the mean curvature as the $L^2$-gradient of the area functional. In particular, $\Sigma$ is a critical point of the area functional if and only if the mean curvature of $\Sigma$ vanishes.

Finally, minimal surfaces can be characterized in terms of the Laplace operator on the submanifold $\Sigma$. To see this, suppose that $a$ is a vector in ambient space $\mathbb{R}^{n+1}$. Let us consider the restriction of the linear function $\langle a,x \rangle$ in $\mathbb{R}^{n+1}$ to the hypersurface $\Sigma$. The Laplacian, on $\Sigma$, of that function is given by 
\[\Delta_\Sigma \langle a,x \rangle = -H \, \langle a,\nu \rangle.\] 
In particular, if $H$ vanishes, then the restriction of the function $\langle a,x \rangle$ to $\Sigma$ gives a harmonic function on $\Sigma$. Conversely, if the restriction of the function $\langle a,x \rangle$ to $\Sigma$ is harmonic for every $a \in \mathbb{R}^{n+1}$, then the mean curvature vanishes identically.

The preceding discussion can be summarized as follows: 

\begin{theorem}
Let $\Sigma$ be a hypersurface in $\mathbb{R}^{n+1}$. Then the following statements are equivalent: 
\begin{itemize}
\item $\Sigma$ is a minimal surface.
\item The first variation of the area functional at $\Sigma$ vanishes.
\item The mean curvature of $\Sigma$ vanishes.
\item Each coordinate function in $\mathbb{R}^{n+1}$ restricts to a harmonic function on $\Sigma$.
\end{itemize}
\end{theorem} 

Two-dimensional minimal surfaces in $\mathbb{R}^3$ can be studied using techniques from complex analysis. In particular, the Enneper--Weierstrass representation makes it possible to describe minimal surfaces locally in terms of holomorphic functions (see \cite{Colding-Minicozzi1}, \cite{Dierkes-Hildebrandt-Sauvigny}). This gives many interesting examples of minimal surfaces in $\mathbb{R}^3$.

While we have so far focused on minimal surfaces in Euclidean space, the notion of a minimal surface makes sense in any ambient Riemannian manifold. The case of minimal surfaces in spheres is particularly interesting: while a minimal surface in Euclidean space can never close up, there do exist examples of closed minimal surfaces in spheres (see \cite{Lawson}). In particular, there is a complete classification of all immersed minimal surfaces in $S^3$ of genus $0$ (cf. \cite{Almgren1}), and of all embedded minimal surfaces in $S^3$ of genus $1$ (cf. \cite{Brendle2}). We refer to \cite{Brendle3} for a survey of some recent developments in this direction. 

Some of the broad themes that have been studied in minimal surface theory are existence and regularity questions; uniqueness questions; and geometric inequalities. In this survey, we will focus on the third topic.

\section{The monotonicity formula}

One of the most fundamental results in minimal surface theory is the monotonicity formula (see e.g. \cite{Simon}). To fix notation, let $B_r$ denote the ball of radius $r$ in the ambient Euclidean space $\mathbb{R}^{n+1}$.

\begin{theorem}
\label{monotonicity}
Let $\Sigma$ be a compact minimal hypersurface in $\mathbb{R}^{n+1}$ with boundary $\partial \Sigma$. Suppose that $\partial \Sigma \cap B_\rho = \emptyset$. Then the function 
\[r \mapsto \frac{|\Sigma \cap B_r|}{|B^n| \, r^n}\] 
is monotone increasing for $0 < r < \rho$.
\end{theorem}

The standard proof of Theorem \ref{monotonicity} uses the co-area formula. In the following, we will present a slightly different argument which is based on an application of the divergence theorem to a suitably chosen vector field. Let us define a vector field $V$ in ambient space $\mathbb{R}^{n+1}$ by 
\[V(x) = |x|^{-n} \, x.\] 
The vector field $V$ has a natural interpretation in terms of the gradient of the Newton potential in $n$-dimensional Euclidean space. For each point $x \in \Sigma$, we denote by $V^{\text{\rm tan}}(x)$ and $V^\perp(x)$ the tangential and normal components of $V(x)$, respectively. The divergence of $V^{\text{\rm tan}}$ is given by 
\begin{align*} 
\text{\rm div}_\Sigma(V^{\text{\rm tan}}) 
&= \sum_{i=1}^n \langle \bar{D}_{e_i} V^{\text{\rm tan}},e_i \rangle \\ 
&= \sum_{i=1}^n \langle \bar{D}_{e_i} V,e_i \rangle - \sum_{i=1}^n \langle \bar{D}_{e_i} V^\perp,e_i \rangle \\ 
&= n \, |x|^{-n} - n \, |x|^{-n-2} \sum_{i=1}^n \langle x,e_i \rangle^2 - H \, \langle V,\nu \rangle \\ 
&= n \, |x|^{-n-2} \, |x^\perp|^2, 
\end{align*} 
where in the last step we have used the fact that $\Sigma$ is minimal. 

By Sard's theorem, there exists a dense open subset $\mathcal{R} \subset (0,\rho)$ such that $\partial B_r$ meets $\Sigma$ transversally for all $r \in \mathcal{R}$. Let us fix two radii $r_0,r_1 \in \mathcal{R}$ such that $r_0 < r_1$. For each point $x \in \Sigma \cap \partial B_{r_1}$, we denote by $\eta$ the co-normal to $\Sigma \cap B_{r_1}$. Similarly, for each point $x \in \Sigma \cap \partial B_{r_0}$, $\eta$ will denote the co-normal to $\Sigma \cap B_{r_0}$. Applying the divergence theorem to the vector field $V^{\text{\rm tan}}$ on $\Sigma \cap (B_{r_1} \setminus B_{r_0})$ gives 
\begin{align*} 
\int_{\Sigma \cap (B_{r_1} \setminus B_{r_0})} n \, |x|^{-n-2} \, |x^\perp|^2 
&= \int_{\Sigma \cap (B_{r_1} \setminus B_{r_0})} \text{\rm div}_\Sigma(V^{\text{\rm tan}}) \\ 
&= \int_{\Sigma \cap \partial B_{r_1}} \langle V,\eta \rangle - \int_{\Sigma \cap \partial B_{r_0}} \langle V,\eta \rangle \\ 
&= r_1^{-n} \int_{\Sigma \cap \partial B_{r_1}} \langle x,\eta \rangle - r_0^{-n} \int_{\Sigma \cap \partial B_{r_0}} \langle x,\eta \rangle. 
\end{align*} 
On the other hand, since $\Sigma$ is minimal, the vector field $x^{\text{\rm tan}}$ satisfies $\text{\rm div}_\Sigma(x^{\text{\rm tan}}) = n$. Applying the divergence theorem to the vector field $x^{\text{\rm tan}}$ gives 
\[n \, |\Sigma \cap B_{r_1}| = \int_{\Sigma \cap B_{r_1}} \text{\rm div}_\Sigma(x^{\text{\rm tan}}) = \int_{\Sigma \cap \partial B_{r_1}} \langle x,\eta \rangle\] 
and 
\[n \, |\Sigma \cap B_{r_0}| = \int_{\Sigma \cap B_{r_0}} \text{\rm div}_\Sigma(x^{\text{\rm tan}}) = \int_{\Sigma \cap \partial B_{r_0}} \langle x,\eta \rangle.\] 
Putting these facts together, we conclude that 
\[\int_{\Sigma \cap (B_{r_1} \setminus B_{r_0})} |x|^{-n-2} \, |x^\perp|^2 = r_1^{-n} \, |\Sigma \cap B_{r_1}| - r_0^{-n} \, |\Sigma \cap B_{r_0}|,\]
which implies the monotonicity formula.

\section{Estimates for the area of a minimal surface in a ball}

One important consequence of Theorem \ref{monotonicity} is that it gives a lower bound for the volume of a minimal surfaces that passes through the origin. 

\begin{theorem}
\label{area.bound.v1}
Let $\Sigma$ be a compact minimal hypersurface in the closed unit ball $\bar{B}^{n+1}$ with boundary $\partial \Sigma \subset \partial B^{n+1}$. If $\Sigma$ passes through the origin, then $|\Sigma| \geq |B^n|$. Moreover, if equality holds, then $\Sigma$ is a flat disk.
\end{theorem}

Theorem \ref{area.bound.v1} follows directly from Theorem \ref{monotonicity}. Indeed, if $\Sigma$ passes through the origin, then $\liminf_{r \to 0} \frac{|\Sigma \cap B_r|}{|B^n| \, r^n} \geq 1$. Hence, the monotonicity formula implies $\frac{|\Sigma \cap B_r|}{|B^n| \, r^n} \geq 1$ for all $0 < r < 1$. From this, the assertion follows. \\

In 1973, Alexander and Osserman \cite{Alexander-Osserman} considered the more general situation when $\Sigma$ passes through some prescribed point in $B^{n+1}$ (not necessarily the origin). They conjectured that, among all minimal surfaces in the unit ball which pass through a prescribed point $y \in B^{n+1}$, the flat disk orthogonal to $y$ has smallest area. This conjecture was proved in \cite{Brendle-Hung}:

\begin{theorem}[S.~Brendle, P.K.~Hung \cite{Brendle-Hung}] 
\label{area.bound.v2}
Let $\Sigma$ be a compact minimal hypersurface in the closed unit ball $\bar{B}^{n+1}$ with boundary $\partial \Sigma \subset \partial B^{n+1}$. If $\Sigma$ passes through a given point $y \in B^{n+1}$, then $|\Sigma| \geq |B^n| \, (1-|y|^2)^{\frac{n}{2}}$. Moreover, if equality holds, then $\Sigma$ is a flat disk which is orthogonal to $y$.
\end{theorem}

Berndtsson \cite{Berndtsson} later gave an alternative proof of Theorem \ref{area.bound.v2} using the theory of supercurrents.

Let us sketch the proof of Theorem \ref{area.bound.v2}. Let us fix a point $y \in B^{n+1}$. We define a vector field $W$ on $\bar{B}^{n+1} \setminus \{y\}$ as follows. For $n>2$, we define  
\begin{align*} 
W(x) 
&= -\frac{1}{n} \, \bigg ( \Big ( \frac{1-2 \langle x,y \rangle+|y|^2}{|x-y|^2} \Big )^{\frac{n}{2}} - 1 \bigg ) \, (x-y) \\ 
&+ \frac{1}{n-2} \, \bigg ( \Big ( \frac{1-2 \langle x,y \rangle+|y|^2}{|x-y|^2} \Big )^{\frac{n-2}{2}} - 1 \bigg ) \, y. 
\end{align*}
For $n=2$, we define  
\begin{align*} 
W(x) 
&= -\frac{1}{2} \, \Big ( \frac{1-2 \langle x,y \rangle+|y|^2}{|x-y|^2} - 1 \Big ) \, (x-y) \\ 
&+ \frac{1}{2} \, \log \Big ( \frac{1-2 \langle x,y \rangle+|y|^2}{|x-y|^2} \Big ) \, y. 
\end{align*}
The vector field $W$ has the following properties: 
\begin{itemize} 
\item For every point $x \in \bar{B}^{n+1}$ and every orthonormal frame $\{e_1,\hdots,e_n\} \subset \mathbb{R}^{n+1}$, we have $\sum_{i=1}^n \langle \bar{D}_{e_i} W,e_i \rangle \leq 1$. 
\item $W(x) = 0$ for $x \in \partial B^{n+1}$.
\item $W(x) = -(1-|y|^2)^{\frac{n}{2}} \, \frac{x-y}{n \, |x-y|^n} + o \big ( \frac{1}{|x-y|^{n-1}} \big )$ as $x \to y$.
\end{itemize}
Since $\Sigma$ is minimal, we have 
\[\text{\rm div}_\Sigma(W^{\text{\rm tan}}) = \sum_{i=1}^n \langle \bar{D}_{e_i} W,e_i \rangle \leq 1,\] 
where $\{e_1,\hdots,e_n\}$ is a local orthonormal frame on $\Sigma$. Applying the divergence theorem to the vector field $W^{\text{\rm tan}}$ on $\Sigma \setminus \{y\}$, we obtain 
\[\int_{\Sigma \setminus \{y\}} \text{\rm div}_\Sigma(W^{\text{\rm tan}}) = |B^n| \, (1-|y|^2)^{\frac{n}{2}}.\] 
Since $\text{\rm div}_\Sigma(W^{\text{\rm tan}}) \leq 1$ at each point on $\Sigma \setminus \{y\}$, we conclude that 
\[|\Sigma| \geq  |B^n| \, (1-|y|^2)^{\frac{n}{2}},\] 
as claimed. \\

Finally, a related argument gives a lower bound for the area of free boundary minimal surfaces in the unit ball: 

\begin{theorem}[S.~Brendle \cite{Brendle1}] 
\label{area.bound.v3}
Let $\Sigma$ be a compact minimal hypersurface in the closed unit ball $\bar{B}^{n+1}$ with boundary $\partial \Sigma \subset \partial B^{n+1}$. If $\Sigma$ meets $\partial B^{n+1}$ orthogonally, then $|\Sigma| \geq |B^n|$. Moreover, if equality holds, then $\Sigma$ is a flat disk.
\end{theorem}

Theorem \ref{area.bound.v3} confirms a conjecture of Schoen. The conjecture has been attributed to a question posed earlier by Guth.

Let us sketch the proof of Theorem \ref{area.bound.v3}; for full details see \cite{Brendle1}. We fix a point $y \in \partial \Sigma$. We define a vector field $Z$ on $\bar{B}^{n+1} \setminus \{y\}$ by 
\[Z(x) = \frac{1}{2} \, x - \frac{x-y}{|x-y|^n} - \frac{n-2}{2} \int_0^1 \frac{tx-y}{|tx-y|^n} \, dt.\] 
The vector field $Z$ has a natural interpretation in terms of the gradient of the Greens function for the Neumann problem on the $n$-dimensional unit ball. The vector field $Z$ has the following properties: 
\begin{itemize}
\item For every point $x \in \bar{B}^{n+1}$ and every orthonormal frame $\{e_1,\hdots,e_n\} \subset \mathbb{R}^{n+1}$, we have $\sum_{i=1}^n \langle \bar{D}_{e_i} Z,e_i \rangle \leq \frac{n}{2}$.
\item The vector field $Z$ is tangential along the boundary $\partial B^{n+1}$; that is, $\langle Z(x),x \rangle = 0$ for all $x \in \partial B^{n+1}$.
\item $Z(x) = -\frac{x-y}{|x-y|^n} + o \big ( \frac{1}{|x-y|^{n-1}} \big )$ as $x \to y$.
\end{itemize}
Since $\Sigma$ is minimal, we have 
\[\text{\rm div}_\Sigma(Z^{\text{\rm tan}}) = \sum_{i=1}^n \langle \bar{D}_{e_i} Z,e_i \rangle \leq \frac{n}{2},\] 
where $\{e_1,\hdots,e_n\}$ is a local orthonormal frame on $\Sigma$. Applying the divergence theorem to the vector field $Z^{\text{\rm tan}}$ on $\Sigma \setminus \{y\}$ gives 
\[\int_{\Sigma \setminus \{y\}} \text{\rm div}_\Sigma(Z^{\text{\rm tan}}) = \frac{n}{2} \, |B^n|.\] 
Since $\text{\rm div}_\Sigma(Z^{\text{\rm tan}}) \leq \frac{n}{2}$ at each point on $\Sigma \setminus \{y\}$, we conclude that 
\[\frac{n}{2} \, |\Sigma| \geq \frac{n}{2} \, |B^n|,\] 
which implies the claim. \\

\section{Gromov's extended monotonicity formula}

The monotonicity formula in Theorem \ref{monotonicity} is valid as long as the ball $B_r$ is disjoint from the boundary $\partial \Sigma$. In this section, we discuss an extended monotonicity formula, due to Gromov \cite{Gromov1}, which holds for all $r$. To fix notation, let $\Sigma$ be a compact minimal hypersurface with boundary $\Gamma = \partial \Sigma$. We assume that the origin does not lie on $\Gamma$. Let $E$ denote the exterior cone over $\Gamma$, i.e. 
\[E = \{\lambda x: x \in \Gamma, \, \lambda \in [1,\infty)\}.\] 
Moreover, we denote by $\tilde{\Sigma} := \Sigma \cup E$ the union of the minimal surface $\Sigma$ and the exterior cone $E$. Note that $\tilde{\Sigma}$ has no boundary, but $\tilde{\Sigma}$ is non-smooth along $\Gamma$. 

\begin{theorem}[M.~Gromov \cite{Gromov1}, Theorem 8.2.A]
\label{extended.monotonicity}
Let $\Sigma$ be a compact minimal hypersurface in $\mathbb{R}^{n+1}$ with boundary $\Gamma = \partial \Sigma$. Suppose that the origin does not lie on $\Gamma$. Let $\tilde{\Sigma}$ denote the extended hypersurface defined above. Then the function 
\[r \mapsto \frac{|\tilde{\Sigma} \cap B_r|}{|B^n| \, r^n}\] 
is monotone increasing for all $r>0$.
\end{theorem}

In the remainder of this section, we sketch the proof of the extended monotonicity formula. As above, we define a vector field $V$ in ambient space $\mathbb{R}^{n+1}$ by 
\[V(x) = |x|^{-n} \, x.\] 
Since $\Sigma$ is minimal, we have $\text{\rm div}_\Sigma(V^{\text{\rm tan}}) = n \, |x|^{-n-2} \, |x^\perp|^2$ and $\text{\rm div}_\Sigma(x^{\text{\rm tan}}) = n$ at each point on $\Sigma$. Moreover, since $E$ is a cone, we have $x^\perp = 0$ at each point on $E$. This implies $\text{\rm div}_E(V^{\text{\rm tan}}) = 0$ and $\text{\rm div}_E(x^{\text{\rm tan}}) = n$ at each point on $E$. 

By Sard's theorem, there exists a dense open subset $\mathcal{R} \subset (0,\infty)$ such that $\partial B_r$ meets $\Sigma$, $E$, and $\Gamma$ transversally for all $r \in \mathcal{R}$. Let us fix two radii $r_0,r_1 \in \mathcal{R}$ such that $r_0 < r_1$. For each point $x \in \Gamma = \partial \Sigma$, we denote by $\eta_\Sigma$ the co-normal to $\Sigma$. Moreover, for each point $x \in \Sigma \cap \partial B_{r_1}$, $\eta_\Sigma$ will denote the co-normal to $\Sigma \cap B_{r_1}$. Finally, for each point $x \in \Sigma \cap \partial B_{r_0}$, $\eta_\Sigma$ will denote the co-normal to $\Sigma \cap B_{r_0}$. Applying the divergence theorem to the vector field $V^{\text{\rm tan}}$ on $\Sigma \cap (B_{r_1} \setminus B_{r_0})$ gives 
\begin{align*} 
\int_{\Sigma \cap (B_{r_1} \setminus B_{r_0})} n \, |x|^{-n-2} \, |x^\perp|^2 
&= r_1^{-n} \int_{\Sigma \cap \partial B_{r_1}} \langle x,\eta_\Sigma \rangle - r_0^{-n} \int_{\Sigma \cap \partial B_{r_0}} \langle x,\eta_\Sigma \rangle \\ 
&+ \int_{\Gamma \cap (B_{r_1} \setminus B_{r_0})} |x|^{-n} \, \langle x,\eta_\Sigma \rangle.
\end{align*} 
On the other hand, applying the divergence theorem to the vector field $x^{\text{\rm tan}}$ on $\Sigma$ gives 
\[n \, |\Sigma \cap B_{r_1}| = \int_{\Sigma \cap \partial B_{r_1}} \langle x,\eta_\Sigma \rangle + \int_{\Gamma \cap B_{r_1}} \langle x,\eta_\Sigma \rangle\] 
and 
\[n \, |\Sigma \cap B_{r_0}| = \int_{\Sigma \cap \partial B_{r_0}} \langle x,\eta_\Sigma \rangle + \int_{\Gamma \cap B_{r_0}} \langle x,\eta_\Sigma \rangle.\] 
Putting these facts together, we obtain 
\begin{align*} 
&\int_{\Sigma \cap (B_{r_1} \setminus B_{r_0})} n \, |x|^{-n-2} \, |x^\perp|^2 \\ 
&= n \, r_1^{-n} \, |\Sigma \cap B_{r_1}| - n \, r_0^{-n} \, |\Sigma \cap B_{r_0}| \\ 
&+ \int_{\Gamma \cap (B_{r_1} \setminus B_{r_0})} (|x|^{-n}-r_1^{-n}) \, \langle x,\eta_\Sigma \rangle + (r_0^{-n}-r_1^{-n}) \int_{\Gamma \cap B_{r_0}} \langle x,\eta_\Sigma \rangle. 
\end{align*} 
We next consider the exterior cone $E$. For each point $x \in \Gamma = \partial E$, we denote by $\eta_E$ the co-normal to $E$. Moreover, for each point $x \in E \cap \partial B_{r_1}$, $\eta_E$ will denote the co-normal to $E \cap B_{r_1}$. Finally, for each point $x \in E \cap \partial B_{r_0}$, $\eta_E$ will denote the co-normal to $E \cap B_{r_0}$. Applying the divergence theorem to the vector field $V^{\text{\rm tan}}$ on $E \cap (B_{r_1} \setminus B_{r_0})$ gives 
\begin{align*} 
0 &= r_1^{-n} \int_{E \cap \partial B_{r_1}} \langle x,\eta_E \rangle - r_0^{-n} \int_{E \cap \partial B_{r_0}} \langle x,\eta_E \rangle \\ 
&+ \int_{\Gamma \cap (B_{r_1} \setminus B_{r_0})} |x|^{-n} \, \langle x,\eta_E \rangle.
\end{align*} 
On the other hand, applying the divergence theorem to the vector field $x^{\text{\rm tan}}$ on $E$ gives 
\[n \, |E \cap B_{r_1}| = \int_{E \cap \partial B_{r_1}} \langle x,\eta_E \rangle + \int_{\Gamma \cap B_{r_1}} \langle x,\eta_E \rangle\] 
and 
\[n \, |E \cap B_{r_0}| = \int_{E \cap \partial B_{r_0}} \langle x,\eta_E \rangle + \int_{\Gamma \cap B_{r_0}} \langle x,\eta_E \rangle.\] 
This implies 
\begin{align*} 
0 &= n \, r_1^{-n} \, |E \cap B_{r_1}| - n \, r_0^{-n} \, |E \cap B_{r_0}| \\ 
&+ \int_{\Gamma \cap (B_{r_1} \setminus B_{r_0})} (|x|^{-n}-r_1^{-n}) \, \langle x,\eta_E \rangle + (r_0^{-n}-r_1^{-n}) \int_{\Gamma \cap B_{r_0}} \langle x,\eta_E \rangle. 
\end{align*} 
We now add the contributions from $\Sigma$ and $E$. This gives 
\begin{align*} 
&\int_{\Sigma \cap (B_{r_1} \setminus B_{r_0})} n \, |x|^{-n-2} \, |x^\perp|^2 \\ 
&= n \, r_1^{-n} \, |\tilde{\Sigma} \cap B_{r_1}| - n \, r_0^{-n} \, |\tilde{\Sigma} \cap B_{r_0}| \\ 
&+ \int_{\Gamma \cap (B_{r_1} \setminus B_{r_0})} (|x|^{-n}-r_1^{-n}) \, \langle x,\eta_\Sigma+\eta_E \rangle + (r_0^{-n}-r_1^{-n}) \int_{\Gamma \cap B_{r_0}} \langle x,\eta_\Sigma+\eta_E \rangle. 
\end{align*} 
We claim that, for each point $x \in \Gamma$, the quantity $\langle x,\eta_\Sigma+\eta_E \rangle$ is nonpositive. Indeed, since $E$ is an exterior cone, the position vector $x$ lies in the tangent space $T_x E$ and is inward-pointing. Consequently, we may write $x = -\lambda \, \eta_E + w$, where $w \in T_x \Gamma$ and $\lambda \geq 0$. Clearly, $\langle w,\eta_\Sigma \rangle = \langle w,\eta_E \rangle = 0$ since $w \in T_x \Gamma$. Hence, we obtain $\langle x,\eta_\Sigma+\eta_E \rangle = -\lambda \, \langle \eta_E,\eta_\Sigma+\eta_E \rangle = -\frac{1}{2} \, \lambda \, |\eta_\Sigma+\eta_E|^2 \leq 0$. Therefore, $\langle x,\eta_\Sigma+\eta_E \rangle$ is nonpositive at each point on $\Gamma$. Thus, 
\[\int_{\Sigma \cap (B_{r_1} \setminus B_{r_0})} |x|^{-n-2} \, |x^\perp|^2 \leq r_1^{-n} \, |\tilde{\Sigma} \cap B_{r_1}| - r_0^{-n} \, |\tilde{\Sigma} \cap B_{r_0}|.\] 
This completes the proof of the extended monotonicitiy formula. \\

In particular, Theorem \ref{extended.monotonicity} implies that the density ratios of $\tilde{\Sigma}$ are bounded from above by the density of $E$ at infinity. In the special case of two-dimensional minimal surfaces in $\mathbb{R}^3$, Ekholm, White, and Wienholtz \cite{Ekholm-White-Wienholtz} were able to estimate the density of $E$ at infinity in terms of the total curvature of $\Gamma$. As a consequence, they obtained the following result:

\begin{theorem}[T.~Ekholm, B.~White, D.~Wienholtz \cite{Ekholm-White-Wienholtz}]
Let $\Sigma$ be a compact minimal surface in $\mathbb{R}^3$ with boundary $\Gamma = \partial \Sigma$. Suppose that the origin does not lie on $\Gamma$. Let $\tilde{\Sigma}$ denote the extended surface defined above. Then 
\[\frac{|\tilde{\Sigma} \cap B_r|}{\pi r^2} \leq \frac{1}{2\pi} \int_\Gamma |k|\] 
for all $r>0$, where $k$ denotes the curvature vector of the boundary $\Gamma$. In particular, if $\int_\Gamma |k| < 4\pi$, then the interior of $\Sigma$ is embedded.
\end{theorem}

\section{The isoperimetric inequality for minimal surfaces} 

\label{abp}

In this section, we discuss sharp isoperimetric inequalities on minimal surfaces. Let us first recall the isoperimetric inequality for domains in Euclidean space: 

\begin{theorem}
\label{isoperimetric.inequality.in.Euclidean.space}
Let $E$ be a compact domain in $\mathbb{R}^n$ with boundary $\partial E$. Then 
\[|\partial E| \geq n \, |B^n|^{\frac{1}{n}} \, |E|^{\frac{n-1}{n}}.\] 
Moreover, if equality holds, then $E$ is a ball.
\end{theorem}

The isoperimetric inequality is one of the most fundamental results in geometry. Many different proofs can be found in the literature. In particular, the isoperimetric inequality is a consequence of the classical Brunn-Minkowski inequality for compact subsets of $\mathbb{R}^n$, which was proved in full generality by Lusternik \cite{Lusternik} in 1935. A modern exposition can be found in \cite{Stein-Shakarchi}, Theorem 5.1. We refer to \cite{Barthe-Cordero-Erausquin} for an alternative proof of the Brunn-Minkowski inequality using heat flows.

The isoperimetric inequality was generalized to the Riemannian setting by Gromov (cf. \cite{Gromov2}, Appendix C). Klartag \cite{Klartag} has developed an alternative approach to the L\'evy-Gromov inequality based on optimal transport and needle decompositions; this approach was generalized to metric measure spaces in \cite{Cavalletti-Mondino}. Moreover, the classical Brunn-Minkowski inequality in Euclidean space is a special case of the Riemannian interpolation inequality proved by Cordero-Erausquin, McCann, and Schmuckenschl\"ager \cite{Cordero-Erausquin-McCann-Schmuckenschlager}. 

It has been conjectured for a long time that the isoperimetric inequality should hold for minimal surfaces. This line of research was initiated in a seminal work of Torsten Carleman in 1921, where he proved a sharp isoperimetric inequality for disk-type minimal surfaces. Various authors have obtained generalizations of this result under weaker topological assumptions (see e.g. \cite{Choe}, \cite{Feinberg}, \cite{Hsiung}, \cite{Li-Schoen-Yau}, \cite{Osserman-Schiffer}, \cite{Reid}). In particular, these results include the case of two-dimensional minimal surfaces with connected boundary:

\begin{theorem}[T.~Carleman \cite{Carleman}; T.~Reid \cite{Reid}; C.C.~Hsiung \cite{Hsiung}]
\label{isoperimetric.inequality.for.2D.minimal.surfaces}
Let $\Sigma$ be a minimal surface in $\mathbb{R}^3$ with boundary $\partial \Sigma$. If $\partial \Sigma$ is connected, then $|\partial \Sigma|^2 \geq 4\pi \, |\Sigma|$. Moreover, if equality holds, then $\Sigma$ is a flat disk.
\end{theorem}

Leon Simon and Andrew Stone have obtained non-sharp isoperimetric inequalities for two-dimensional minimal surfaces (see \cite{Stone} and \cite{Topping}, Section 4). These results require no topological assumptions.

In the following, we present the proof of Theorem \ref{isoperimetric.inequality.for.2D.minimal.surfaces}. The proof is a generalization of Hurwitz's proof of the isoperimetric inequality in $\mathbb{R}^2$ (see \cite{Hurwitz}, pp.~392--394). By scaling, we may assume that $|\partial \Sigma| = 2\pi$. By assumption, $\partial \Sigma$ is connected. Let $\alpha: [0,2\pi] \to \partial \Sigma$ denote a parametrization of $\partial \Sigma$ by arclength, so that $|\alpha'(s)|=1$ for all $s \in [0,2\pi]$. Without loss of generality, we may assume that the center of mass of the boundary $\partial \Sigma$ is at the origin, so that $\int_0^{2\pi} \alpha_i(s) \, ds = 0$ for each $1 \leq i \leq 3$. Applying Wirtinger's inequality to the function $\alpha_i(s)$, we obtain 
\[\int_0^{2\pi} \alpha_i(s)^2 \, ds \leq \int_0^{2\pi} \alpha_i'(s)^2 \, ds\] 
for each $1 \leq i \leq 3$. Summation over $i$ gives 
\[\int_0^{2\pi} |\alpha(s)|^2 \, ds \leq \int_0^{2\pi} |\alpha'(s)|^2 \, ds = 2\pi.\] 
In other words, 
\[\int_{\partial \Sigma} |x|^2 \leq 2\pi.\] 
On the other hand, since $\Sigma$ is minimal, it follows that $\text{\rm div}_\Sigma(x^{\text{\rm tan}}) = 2$. Hence, the divergence theorem gives 
\[4 \, |\Sigma| = 2 \int_\Sigma \text{\rm div}_\Sigma(x^{\text{\rm tan}}) = 2 \int_{\partial \Sigma} \langle x,\eta \rangle \leq \int_{\partial \Sigma} |x|^2 + \int_{\partial \Sigma} |\eta|^2 \leq 4\pi.\] 
Here, $\eta$ denotes the co-normal to $\Sigma$; in particular, $|\eta|=1$ at each point on $\partial \Sigma$. Thus, $|\Sigma| \leq \pi$, which implies the isoperimetric inequality. 

Finally, we give the proof of the rigidity statement. Suppose that equality holds in the isoperimetric inequality. By scaling, we can arrange that $|\partial \Sigma| = 2\pi$ and $|\Sigma| = \pi$. As above, we assume that the center of mass of the boundary $\partial \Sigma$ is at the origin, and that $\alpha: [0,2\pi] \to \partial \Sigma$ is a parametrization of $\partial \Sigma$ by arclength. For each $i$, the function $\alpha_i(s)$ must achieve equality in Wirtinger's inequality. This implies 
\[\alpha(s) = \cos(s) \, v + \sin(s) \, w\] 
for all $s \in [0,2\pi]$, where $v$ and $w$ are fixed vectors in $\mathbb{R}^3$. This gives 
\begin{align*} 
1 
&= |\alpha'(s)|^2 \\ 
&= \sin^2(s) \, |v|^2 + \cos^2(s) \, |w|^2 - 2 \sin(s) \cos(s) \, \langle v,w \rangle \\ 
&= \frac{1}{2} \, (|v|^2+|w|^2) - \frac{1}{2} \, \cos(2s) \, (|v|^2-|w|^2) - \sin (2s) \, \langle v,w \rangle 
\end{align*}
for all $s \in [0,2\pi]$. Consequently, $|v|^2+|w|^2=2$, $|v|^2-|w|^2 = 0$, and $\langle v,w \rangle = 0$. Therefore, $v$ and $w$ are orthonormal, and $\partial \Sigma$ is a circle of radius $1$ which lies in the plane spanned by $v$ and $w$. This completes the proof of Theorem \ref{isoperimetric.inequality.for.2D.minimal.surfaces}. \\

We now turn to the higher dimensional case. A fundamental result in higher dimensions is the Michael-Simon Sobolev inequality (cf. \cite{Allard}, Section 7, and \cite{Michael-Simon}). This inequality holds for an arbitrary hypersurface in Euclidean space. It implies an isoperimetric inequality for minimal surfaces, albeit with a non-sharp constant. Castillon \cite{Castillon} later gave an alternative proof of the Michael-Simon Sobolev inequality using ideas from optimal transport; again, this gives a non-sharp constant. In a recent paper \cite{Brendle4}, we obtained a sharp version of the Michael-Simon Sobolev inequality.

\begin{theorem}[S.~Brendle \cite{Brendle4}] 
\label{sharp.version.of.michael.simon}
Let $\Sigma$ be a compact hypersurface in $\mathbb{R}^{n+1}$ (possibly with boundary $\partial \Sigma$), and let $f$ be a positive smooth function on $\Sigma$. Then 
\[\int_\Sigma \sqrt{|\nabla^\Sigma f|^2 + f^2 H^2} + \int_{\partial \Sigma} f \geq n \, |B^n|^{\frac{1}{n}} \, \Big ( \int_\Sigma f^{\frac{n}{n-1}} \Big )^{\frac{n-1}{n}}.\] 
Moreover, if equality holds, then $f$ is constant and $\Sigma$ is a flat disk.
\end{theorem}

Theorem \ref{sharp.version.of.michael.simon} actually holds for every submanifold of codimension at most $2$. If the codimension is $3$ or higher, only a non-sharp version of the inequality is known. A similar inequality holds in Riemannian manifolds with nonnegative sectional curvature (cf. \cite{Brendle6}); in that case, the constant in the inequality depends not only on the dimension, but also on the asymptotic volume ratio of the ambient manifold.

The proof of Theorem \ref{sharp.version.of.michael.simon} uses the Alexandrov-Bakelman-Pucci method and is inspired in part by an elegant argument due to Cabr\'e \cite{Cabre1},\cite{Cabre2} (see also \cite{Trudinger}). In the following, we describe the main ideas in the codimension $1$ case; we refer to \cite{Brendle4} for a detailed proof in the codimension $2$ setting. First, it is enough to prove the assertion in the special case when $\Sigma$ is connected. (If $\Sigma$ is disconnected, we apply the inequality to each connected component, and take the sum over all connected components.) Second, by scaling, it is enough to prove the assertion in the special case when 
\[\int_\Sigma \sqrt{|\nabla^\Sigma f|^2 + f^2 H^2} + \int_{\partial \Sigma} f = n \int_\Sigma f^{\frac{n}{n-1}}.\] 
This normalization ensures that we can find a function $u: \Sigma \to \mathbb{R}$ which solves the PDE 
\[\text{\rm div}_\Sigma(f \, \nabla^\Sigma u) = n \, f^{\frac{n}{n-1}} - \sqrt{|\nabla^\Sigma f|^2 + f^2 H^2}\] 
on $\Sigma$ with Neumann boundary condition $\langle \nabla^\Sigma u,\eta \rangle = 1$ on $\partial \Sigma$. Here, $\eta$ denotes the co-normal to $\Sigma$. Note that $u$ is of class $C^{2,\gamma}$ for each $0 < \gamma < 1$ by standard elliptic regularity theory. 

Let
\begin{align*} 
\Omega &:= \{x \in \Sigma \setminus \partial \Sigma: |\nabla^\Sigma u(x)| < 1\}, \\ 
U &:= \{(x,y): x \in \Sigma \setminus \partial \Sigma, \, y \in T_x^\perp \Sigma, \, |\nabla^\Sigma u(x)|^2 + |y|^2 < 1\}, \\ 
A &:= \{(x,y) \in U: D_\Sigma^2 u(x) + h(x) \, \langle \nu(x),y \rangle \geq 0\}, 
\end{align*} 
where $h$ denotes the second fundamental form of $\Sigma$. We define a map $\Phi: U \to \mathbb{R}^{n+1}$ by 
\[\Phi(x,y) = \nabla^\Sigma u(x) + y\] 
for all $(x,y) \in U$. One can show that the image $\Phi(A)$ is the open unit ball $B^{n+1}$. The Jacobian determinant of $\Phi$ satisfies 
\[\det D\Phi(x,y) = \det (D_\Sigma^2 u(x) + h(x) \, \langle \nu(x),y \rangle)\] 
for all $(x,y) \in U$. Using the PDE for $u$, we obtain 
\begin{align*} 
\Delta_\Sigma u(x) 
&= n \, f(x)^{\frac{1}{n-1}} - f(x)^{-1} \, \langle \nabla^\Sigma f(x),\nabla^\Sigma u(x) \rangle \\ 
&- f(x)^{-1} \, \sqrt{|\nabla^\Sigma f(x)|^2 + f(x)^2 H(x)^2} 
\end{align*}
for all $x \in \Sigma$. Since $|\nabla^\Sigma u(x)|^2+|y|^2 < 1$ for all $(x,y) \in U$, the Cauchy-Schwarz inequality gives  
\[-\langle \nabla^\Sigma f(x),\nabla^\Sigma u(x) \rangle + f(x) \, H(x) \, \langle \nu(x),y \rangle \leq \sqrt{|\nabla^\Sigma f(x)|^2 + f(x)^2 H(x)^2}\] 
for all $(x,y) \in U$. Putting these facts together, we obtain 
\[\Delta_\Sigma u(x) + H(x) \, \langle \nu(x),y \rangle \leq n \, f(x)^{\frac{1}{n-1}}\] 
for all $(x,y) \in U$. Since $D_\Sigma^2 u(x) + h(x) \, \langle \nu(x),y \rangle \geq 0$ for all $(x,y) \in A$, the arithmetic-geometric mean inequality implies
\[0 \leq \det (D_\Sigma^2 u(x) + h(x) \, \langle \nu(x),y \rangle) \leq \Big ( \frac{\Delta_\Sigma u(x) + H(x) \, \langle \nu(x),y \rangle}{n} \Big )^n \leq f(x)^{\frac{n}{n-1}}\] 
for all $(x,y) \in A$. Therefore, 
\[0 \leq \det D\Phi(x,y) \leq f(x)^{\frac{n}{n-1}}\] 
for all $(x,y) \in A$. We now apply the change of variables formula to the map $\Phi$. This gives 
\begin{align*} 
\pi \, |B^n| 
&= \int_{B^{n+1}} \frac{1}{\sqrt{1-|\xi|^2}} \, d\xi \\ 
&\leq \int_\Omega \bigg ( \int_{\{y \in T_x^\perp \Sigma: |\Phi(x,y)|^2 < 1\}} \frac{|\det D\Phi(x,y)|}{\sqrt{1-|\Phi(x,y)|^2}} \, 1_A(x,y) \, dy \bigg ) \, d\text{\rm vol}(x) \\ 
&\leq \int_\Omega \bigg ( \int_{\{y \in T_x^\perp \Sigma: |\nabla^\Sigma u(x)|^2+|y|^2 < 1\}} \frac{f(x)^{\frac{n}{n-1}}}{\sqrt{1-|\nabla^\Sigma u(x)|^2-|y|^2}} \, dy \bigg ) \, d\text{\rm vol}(x) \\ 
&= \pi \int_\Omega f(x)^{\frac{n}{n-1}} \, d\text{\rm vol}(x).
\end{align*} 
In the last step, we have used the fact that the normal space $T_x^\perp \Sigma$ is one-dimensional and $\int_{-a}^a \frac{1}{\sqrt{a^2-y^2}} \, dy = \pi$ for each $a > 0$. Consequently, 
\[|B^n| \leq \int_\Omega f^{\frac{n}{n-1}} \leq \int_\Sigma f^{\frac{n}{n-1}}.\] 
Thus, we conclude that 
\[\int_\Sigma \sqrt{|\nabla^\Sigma f|^2 + f^2 H^2} + \int_{\partial \Sigma} f = n \int_\Sigma f^{\frac{n}{n-1}} \geq n \, |B^n|^{\frac{1}{n}} \, \Big ( \int_\Sigma f^{\frac{n}{n-1}} \Big )^{\frac{n-1}{n}},\]
as claimed. 

Finally, we sketch the proof of the rigidity statement. Suppose that equality holds. It is easy to see that $\Sigma$ must be connected. By scaling, we can arrange that $\int_\Sigma \sqrt{|\nabla^\Sigma f|^2 + f^2 H^2} + \int_{\partial \Sigma} f = n \, |B^n|$ and $\int_\Sigma f^{\frac{n}{n-1}} = |B^n|$. Let $u: \Sigma \to \mathbb{R}$, $\Omega$, $U$, $A$, and $\Phi: U \to \mathbb{R}^{n+1}$ be defined as above. We first observe that $\int_\Omega  f^{\frac{n}{n-1}} = \int_\Sigma f^{\frac{n}{n-1}}$. Consequently, the complement $\Sigma \setminus \Omega$ has $n$-dimensional measure zero. Moreover, the set 
\[U \setminus \{(x,y) \in A: \det D\Phi(x,y) = f(x)^{\frac{n}{n-1}}\}\] 
has $(n+1)$-dimensional measure zero. Now, if $(x,y) \in A$ and $\det D\Phi(x,y) = f(x)^{\frac{n}{n-1}}$, then equality holds in the arithmetic-geometric mean inequality, and this implies $D_\Sigma^2 u(x) + h(x) \, \langle \nu(x),y \rangle = f(x)^{\frac{1}{n-1}} \, g$. Consequently, the set 
\[U \setminus \{(x,y) \in A: D_\Sigma^2 u(x) + h(x) \, \langle \nu(x),y \rangle = f(x)^{\frac{1}{n-1}} \, g\}\] 
has $(n+1)$-dimensional measure zero. Therefore, $D_\Sigma^2 u(x) + h(x) \, \langle \nu(x),y \rangle = f(x)^{\frac{1}{n-1}} \, g$ for all points $(x,y) \in U$. This implies $D_\Sigma^2 u = f^{\frac{1}{n-1}} \, g$ and $h = 0$ at each point in $\Omega$. Using the PDE for $u$, we obtain $\langle \nabla^\Sigma f,\nabla^\Sigma u \rangle = -|\nabla^\Sigma f|$, hence $\nabla^\Sigma f = 0$ at each point in $\Omega$. Since $\Omega$ is a dense subset of $\Sigma$, we conclude that $D_\Sigma^2 u = f^{\frac{1}{n-1}} \, g$, $h = 0$, $\nabla^\Sigma f = 0$, and $|\nabla^\Sigma u| \leq 1$ at each point in $\Sigma$. To summarize, $\Sigma$ is contained in a hyperplane $P$; we have $f=\lambda^{n-1}$ for some positive constant $\lambda$; and the function $u$ is of the form $u(x) = \frac{1}{2} \, \lambda \, |x-p|^2 + c$ for some point $p \in P$ and some constant $c$. Since $|\nabla^\Sigma u| \leq 1$ at each point on $\Sigma$, it follows that $\Sigma$ is contained in the intersection of $P$ with a closed ball of radius $\lambda^{-1}$ around $p$. On the other hand, since $\int_\Sigma f^{\frac{n}{n-1}} = |B^n|$, the volume of $\Sigma$ is given by $|B^n| \, \lambda^{-n}$. Thus, $\Sigma$ is the intersection of $P$ with a closed ball of radius $\lambda^{-1}$. This completes the proof of the rigidity statement. \\

\begin{corollary}[S.~Brendle \cite{Brendle4}]
\label{isoperimetric.inequality}
Let $\Sigma$ be a compact minimal hypersurface in $\mathbb{R}^{n+1}$ with boundary $\partial \Sigma$. Then 
\[|\partial \Sigma| \geq n \, |B^n|^{\frac{1}{n}} \, |\Sigma|^{\frac{n-1}{n}}.\] 
Moreover, if equality holds, then $\Sigma$ is a flat disk.
\end{corollary}

Almgren \cite{Almgren2} has obtained a sharp version of the filling inequality of Federer and Fleming \cite{Federer-Fleming}. As a consequence, he was able to prove the sharp isoperimetric inequality under the assumption that $\Sigma$ is an absolute minimizer of area.

In the remainder of this section, we discuss several consequences of Corollary \ref{isoperimetric.inequality}.

First, Corollary \ref{isoperimetric.inequality} gives a lower bound for the area of a minimal surface in the unit ball under an assumption on the contact angle.

\begin{theorem}
\label{area.bound.v4}
Fix a real number $\theta \in (0,\frac{\pi}{2}]$. Let $\Sigma$ be a compact minimal hypersurface in the closed unit ball $\bar{B}^{n+1}$ with boundary $\partial \Sigma \subset \partial B^{n+1}$. Suppose that, at each point on $\partial \Sigma$, the contact angle between $\Sigma$ and $\partial B^{n+1}$ is at least $\theta$, so that $|x^\perp| \leq \cos \theta$ for all $x \in \partial \Sigma$. Then 
\[|\Sigma| \geq |B^n| \, \sin^n \theta.\] 
Moreover, if equality holds, then $\Sigma$ is a flat disk.
\end{theorem}

If $\Sigma$ meets $\partial B^{n+1}$ orthogonally, then we may apply Theorem \ref{area.bound.v4} with $\theta = \frac{\pi}{2}$. In this special case, the statement of Theorem \ref{area.bound.v4} reduces to Theorem \ref{area.bound.v3}. 

Let us indicate how Theorem \ref{area.bound.v4} follows from Corollary \ref{isoperimetric.inequality}. For each point $x \in \partial \Sigma$, we denote by $\eta$ the co-normal to $\Sigma$. Clearly, $\eta = \frac{x^{\text{\rm tan}}}{|x^{\text{\rm tan}}|}$. The assumption on the contact angle implies $|x^\perp| \leq \cos \theta$ for all $x \in \partial \Sigma$. Moreover, $|x^\perp|^2 + |x^{\text{\rm tan}}|^2 = |x|^2 = 1$ for all $x \in \partial \Sigma$. Consequently, $|x^{\text{\rm tan}}| \geq \sin \theta$ for all $x \in \partial \Sigma$. This implies $\langle x,\eta \rangle = |x^{\text{\rm tan}}| \geq \sin \theta$ for all $x \in \partial \Sigma$. Using the formula $\text{\rm div}_\Sigma(x^{\text{\rm tan}}) = n$ and the divergence theorem, we obtain 
\[n \, |\Sigma| = \int_\Sigma \text{\rm div}_\Sigma(x^{\text{\rm tan}}) = \int_{\partial \Sigma} \langle x,\eta \rangle \geq |\partial \Sigma| \, \sin \theta.\] 
Moreover, Corollary \ref{isoperimetric.inequality} implies 
\[|\partial \Sigma| \geq n \, |B^n|^{\frac{1}{n}} \, |\Sigma|^{\frac{n-1}{n}}.\] 
Putting these facts together, we conclude that 
\[n \, |\Sigma| \geq n \, |B^n|^{\frac{1}{n}} \, |\Sigma|^{\frac{n-1}{n}} \, \sin \theta,\] 
which implies the claim. \\

Second, using Corollary \ref{isoperimetric.inequality} we obtain a Brunn-Minkowski-type inequality on minimal hypersurfaces. Recall that the classical Brunn-Minkowski inequality gives a lower bound for the volume of a tubular neighborhood of a compact subset of $\mathbb{R}^n$. More precisely, if $E$ is a compact subset of $\mathbb{R}^n$, then the volume of the tubular neighborhood $E_r := \{x \in \mathbb{R}^n: \inf_{y \in E} |x-y| \leq r\}$ can be estimated by $|E_r|^{\frac{1}{n}} \geq |E|^{\frac{1}{n}} + |B^n|^{\frac{1}{n}} \, r$. In the following, we extend this inequality to the setting of minimal hypersurfaces.

\begin{theorem}
\label{brunn.minkowski.intrinsic.version}
Let $\Sigma$ be a compact minimal hypersurface in $\mathbb{R}^{n+1}$ with boundary $\partial \Sigma$. Let $E$ be a compact subset of $\Sigma$, and let $\mathcal{E}_r$ denote the set of all points in $\Sigma$ which have intrinsic distance at most $r$ from the set $E$. Moreover, suppose that the intrinsic distance of the set $E$ from the boundary $\partial \Sigma$ is greater than $\rho$. Then 
\[|\mathcal{E}_r|^{\frac{1}{n}} \geq |E|^{\frac{1}{n}} + |B^n|^{\frac{1}{n}} \, r\] 
for $0 < r < \rho$. 
\end{theorem} 

Theorem \ref{brunn.minkowski.intrinsic.version} follows by combining Corollary \ref{isoperimetric.inequality} with the co-area formula. To explain this, let us fix a radius $r$ such that $0 < r < \rho$, and let $f: \Sigma \to \mathbb{R}$ denote the intrinsic distance from the set $E$. Clearly, $f$ is Lipschitz continuous with Lipschitz constant $1$. Using the convolution technique of Greene and Wu (see \cite{Greene-Wu}, Section 2), we can construct a sequence of smooth functions $f_j$ with the following properties: 
\begin{itemize} 
\item The function $f_j$ is defined on an open subset $\Omega_j$ of $\Sigma$, and $\Omega_j$ contains the set $\mathcal{E}_\rho$. 
\item $\sup_{\Omega_j} |\nabla^\Sigma f_j| \leq 1+\delta_j$, where $\delta_j \to 0$. 
\item $\sup_{\Omega_j} |f_j-f| \leq \varepsilon_j$, where $\varepsilon_j \to 0$. 
\end{itemize} 
In the following, we choose $j$ sufficiently large so that $2\varepsilon_j < r$. Since $|\nabla^\Sigma f_j| \leq 1+\delta_j$, the co-area formula gives 
\[\frac{d}{ds} |\{x \in \Omega_j: f_j(x) \leq s\}| \geq (1+\delta_j)^{-1} \, |\{x \in \Omega_j: f_j(x)=s\}|\] 
whenever $s \in (\varepsilon_j,r-\varepsilon_j)$ is a regular value of $f_j$. Moreover, Corollary \ref{isoperimetric.inequality} implies 
\[|\{x \in \Omega_j: f_j(x)=s\}| \geq n \, |B^n|^{\frac{1}{n}} \, |\{x \in \Omega_j: f_j(x) \leq s\}|^{\frac{n-1}{n}}\] 
whenever $s \in (\varepsilon_j,r-\varepsilon_j)$ is a regular value of $f_j$. Putting these facts together, we obtain 
\[\frac{d}{ds} |\{x \in \Omega_j: f_j(x) \leq s\}| \geq (1+\delta_j)^{-1} \, n \, |B^n|^{\frac{1}{n}} \, |\{x \in \Omega_j: f_j(x) \leq s\}|^{\frac{n-1}{n}},\] 
hence 
\[\frac{d}{ds} |\{x \in \Omega_j: f_j(x) \leq s\}|^{\frac{1}{n}} \geq (1+\delta_j)^{-1} \, |B^n|^{\frac{1}{n}}\] 
whenever $s \in (\varepsilon_j,r-\varepsilon_j)$ is a regular value of $f_j$. Since the function $s \mapsto |\{x \in \Omega_j: f_j(x) \leq s\}|^{\frac{1}{n}}$ is monotone increasing, we conclude that 
\begin{align*} 
&|\{x \in \Omega_j: f_j(x) \leq r-\varepsilon_j\}|^{\frac{1}{n}} -  |\{x \in \Omega_j: f_j(x) \leq \varepsilon_j\}|^{\frac{1}{n}} \\ 
&\geq \int_{\varepsilon_j}^{r-\varepsilon_j} \frac{d}{ds} |\{x \in \Omega_j: f_j(x) \leq s\}|^{\frac{1}{n}} \, ds \\ 
&\geq (1+\delta_j)^{-1} \, |B^n|^{\frac{1}{n}} \, (r-2\varepsilon_j). 
\end{align*} 
We next observe that $E \subset \{x \in \Omega_j: f(x) = 0\} \subset \{x \in \Omega_j: f_j(x) \leq \varepsilon_j\}$ and $\{x \in \Omega_j: f_j(x) \leq r-\varepsilon_j\} \subset \{x \in \Omega_j: f(x) \leq r\} \subset \mathcal{E}_r$. Consequently, 
\[|\mathcal{E}_r|^{\frac{1}{n}} -  |E|^{\frac{1}{n}} \geq (1+\delta_j)^{-1} \, |B^n|^{\frac{1}{n}} \, (r-2\varepsilon_j).\] 
Passing to the limit as $j \to \infty$ gives 
\[|\mathcal{E}_r|^{\frac{1}{n}} - |E|^{\frac{1}{n}} \geq |B^n|^{\frac{1}{n}} \, r,\] 
as claimed. \\

\begin{corollary}
\label{brunn.minkowski.extrinsic.version}
Let $\Sigma$ be a compact minimal hypersurface in $\mathbb{R}^{n+1}$ with boundary $\partial \Sigma$. Let $E$ be a compact subset of $\Sigma$, and let $E_r := \{x \in \mathbb{R}^{n+1}: \inf_{y \in E} |x-y| \leq r\}$ denote the set of all points in ambient space $\mathbb{R}^{n+1}$ which have distance at most $r$ from the set $E$. Moreover, suppose that $\partial \Sigma \cap E_\rho = \emptyset$. Then 
\[|\Sigma \cap E_r|^{\frac{1}{n}} \geq |E|^{\frac{1}{n}} + |B^n|^{\frac{1}{n}} \, r\] 
for $0 < r < \rho$. 
\end{corollary} 

Since $\mathcal{E}_r \subset \Sigma \cap E_r$, Corollary \ref{brunn.minkowski.extrinsic.version} is a direct consequence of Theorem \ref{brunn.minkowski.intrinsic.version}. 

In the special case when $E$ consists of a single point, Corollary \ref{brunn.minkowski.extrinsic.version} gives an alternative proof of Theorem \ref{area.bound.v1}.

Finally, Corollary \ref{isoperimetric.inequality} implies that the sharp $L^p$ Sobolev inequality of Aubin \cite{Aubin} and Talenti \cite{Talenti} holds on every minimal hypersurface:

\begin{theorem}
\label{Lp.Sobolev.inequality}
Let $\Sigma$ be a compact minimal hypersurface in $\mathbb{R}^{n+1}$ with boundary $\partial \Sigma$, and let $1 < p < n$. Let $f$ be a nonnegative smooth function on $\Sigma$ which vanishes in a neighborhood of $\partial \Sigma$. Then 
\[\bigg ( \int_\Sigma f^{\frac{np}{n-p}} \bigg )^{\frac{n-p}{np}} \leq K(n,p) \, \bigg ( \int_\Sigma |\nabla^\Sigma f|^p \bigg )^{\frac{1}{p}},\] 
where 
\[K(n,p) = \pi^{-\frac{1}{2}} \, n^{-\frac{1}{p}} \, \Big ( \frac{p-1}{n-p} \Big )^{\frac{p-1}{p}} \, \Big ( \frac{\Gamma(\frac{n}{2}+1) \, \Gamma(n)}{\Gamma(\frac{n}{p}) \, \Gamma(n+1-\frac{n}{p})} \Big )^{\frac{1}{n}}.\] 
\end{theorem}

Theorem \ref{Lp.Sobolev.inequality} follows by combining the isoperimetric inequality in Corollary \ref{isoperimetric.inequality} with the co-area formula. The argument is the same as in \cite{Aubin} and \cite{Talenti}.

\section{The logarithmic Sobolev inequality on a self-similar shrinker}

The classical logarithmic Sobolev inequality in Euclidean space has been studied by many authors (see e.g. \cite{Bakry}, \cite{Bakry-Emery}, \cite{Bobkov-Ledoux}, \cite{Cordero-Erausquin}, \cite{Gross1}, \cite{Gross2}, \cite{Ledoux1}, \cite{Ledoux2}). The statement is as follows: 

\begin{theorem}[L.~Gross \cite{Gross1}]
\label{log.Sobolev.inequality.in.Gauss.space}
Let 
\[d\gamma = (4\pi)^{-\frac{n}{2}} \, e^{-\frac{|x|^2}{4}} \, dx\] 
denote the Gaussian measure on $\mathbb{R}^n$. Then 
\[\int_{\mathbb{R}^n} \varphi \, \log \varphi \, d\gamma - \int_{\mathbb{R}^n} \frac{|\nabla \varphi|^2}{\varphi} \, d\gamma \leq \bigg ( \int_{\mathbb{R}^n} \varphi \, d\gamma \bigg ) \, \log \bigg ( \int_{\mathbb{R}^n} \varphi \, d\gamma \bigg )\] 
for every positive smooth function $\varphi$ on $\mathbb{R}^n$ satisfying $\int_{\mathbb{R}^n} \varphi \, d\gamma < \infty$ and $\int_{\mathbb{R}^n} \frac{|\nabla \varphi|^2}{\varphi} \, d\gamma < \infty$.
\end{theorem}

There are many different proofs of the logarithmic Sobolev inequality. These employ a variety of techniques, including the central limit theorem \cite{Gross2}, heat flows \cite{Bakry-Emery}, and optimal transport \cite{Cordero-Erausquin}. The logarithmic Sobolev inequality can be viewed as a corollary of the isoperimetric inequality in Gauss space (see \cite{Borell1}, \cite{Borell2}, \cite{Sudakov-Tsirelson}). The Gaussian isoperimetric inequality can be proven using the heat equation; see \cite{Ledoux2}, Section 1.2. 

The logarithmic Sobolev inequality in Euclidean space is a special case of the Bakry-\'Emery theorem \cite{Bakry-Emery}. Similarly, the isoperimetric inequality in Gauss space is a special case of the isoperimetric comparison theorem of Bakry-Ledoux \cite{Bakry-Ledoux}.

Ecker proved a logarithmic Sobolev inequality which holds on every submanifold of Euclidean space, albeit with a non-sharp constant (see \cite{Ecker1} and \cite{Ecker2}, pp.~59--60). This inequality is similar in spirit to the Michael-Simon Sobolev inequality. Using the techniques in Section \ref{abp}, we obtain a sharp version of Ecker's inequality:

\begin{theorem}[S.~Brendle \cite{Brendle5}]
\label{log.Sobolev}
Let $\Sigma$ be a compact hypersurface in $\mathbb{R}^{n+1}$ without boundary, and let 
\[d\gamma = (4\pi)^{-\frac{n}{2}} \, e^{-\frac{|x|^2}{4}} \, d\text{\rm vol}\] 
denote the Gaussian measure on $\Sigma$. Then 
\begin{align*} 
&\int_\Sigma \varphi \, \log \varphi \, d\gamma - \int_\Sigma \frac{|\nabla^\Sigma \varphi|^2}{\varphi} \, d\gamma - \int_\Sigma \varphi \, \Big ( H-\frac{1}{2} \, \langle x,\nu \rangle \Big )^2 \, d\gamma \\ 
&\leq \bigg ( \int_\Sigma \varphi \, d\gamma \bigg ) \, \log \bigg ( \int_\Sigma \varphi \, d\gamma \bigg ) 
\end{align*}
for every positive smooth function $\varphi$ on $\Sigma$.
\end{theorem}

Theorem \ref{log.Sobolev} actually holds for submanifolds of arbitrary codimension; see \cite{Brendle5}.

Let us sketch the proof of Theorem \ref{log.Sobolev}. As in Section \ref{abp}, we can reduce to the special case when $\Sigma$ is connected. By scaling, we may assume that 
\[\int_\Sigma \varphi \, \log \varphi \, d\gamma - \int_\Sigma \frac{|\nabla^\Sigma \varphi|^2}{\varphi} \, d\gamma - \int_\Sigma \varphi \, \Big ( H-\frac{1}{2} \, \langle x,\nu \rangle \Big )^2 \, d\gamma = 0.\] 
This normalization ensures that we can find a smooth function $v: \Sigma \to \mathbb{R}$ such that 
\begin{align*} 
&\text{\rm div}_\Sigma(e^{-\frac{|x|^2}{4}} \, \varphi \, \nabla^\Sigma v) \\ 
&= e^{-\frac{|x|^2}{4}} \, \varphi \, \log \varphi - e^{-\frac{|x|^2}{4}} \, \frac{|\nabla^\Sigma \varphi|^2}{\varphi} - e^{-\frac{|x|^2}{4}} \, \varphi \, \Big ( H-\frac{1}{2} \, \langle x,\nu \rangle \Big )^2. 
\end{align*} 
Let $u(x) := v(x) + \frac{|x|^2}{2}$. We define 
\begin{align*} 
U &:= \{(x,y): x \in \Sigma, \, y \in T_x^\perp \Sigma\}, \\ 
A &:= \{(x,y) \in U: D_\Sigma^2 u(x) + h(x) \, \langle \nu(x),y \rangle \geq 0\}, 
\end{align*} 
where $h$ denotes the second fundamental form of $\Sigma$. Moreover, we define a map $\Phi: U \to \mathbb{R}^{n+1}$ by 
\[\Phi(x,y) = \nabla^\Sigma u(x) + y\] 
for all $(x,y) \in U$. It can be shown that the image $\Phi(A)$ is all of $\mathbb{R}^{n+1}$. The Jacobian determinant of $\Phi$ satisfies 
\[\det D\Phi(x,y) = \det (D_\Sigma^2 u(x) + h(x) \, \langle \nu(x),y \rangle)\] 
for all $(x,y) \in U$. Using the PDE for $v$, we obtain 
\begin{align*}  
&\Delta_\Sigma u(x) - \frac{|\nabla^\Sigma u(x)|^2}{4} + \frac{|x|^2}{4} + H(x)^2 - n \\ 
&= \log \varphi(x) - \frac{|2 \, \nabla^\Sigma \varphi(x) + \varphi(x) \, \nabla^\Sigma v(x)|^2}{4\varphi(x)^2} \leq \log \varphi(x) 
\end{align*} 
for all $x \in \Sigma$. This implies  
\begin{align*} 
&\Delta_\Sigma u(x) + H(x) \, \langle \nu(x),y \rangle - n \\ 
&\leq \frac{|\nabla^\Sigma u(x)|^2+|y|^2}{4} -\frac{|x|^2}{4} - \frac{|y-2H(x)\nu(x)|^2}{4} + \log \varphi(x) 
\end{align*} 
for all $(x,y) \in U$. Since $D_\Sigma^2 u(x) + h(x) \, \langle \nu(x),y \rangle \geq 0$ for all $(x,y) \in A$, it follows that 
\begin{align*} 
0 
&\leq e^{-\frac{|\nabla^\Sigma u(x)|^2+|y|^2}{4}} \, \det (D_\Sigma^2 u(x) + h(x) \, \langle \nu(x),y \rangle) \\ 
&\leq e^{-\frac{|\nabla^\Sigma u(x)|^2+|y|^2}{4}} \, e^{\Delta_\Sigma u(x) + H(x) \, \langle \nu(x),y \rangle - n} \\ 
&\leq e^{-\frac{|x|^2}{4} - \frac{|y-2H(x)\nu(x)|^2}{4}} \, \varphi(x) 
\end{align*}
for all $(x,y) \in A$. Thus, we conclude that 
\[0 \leq e^{-\frac{|\Phi(x,y)|^2}{4}} \, \det D\Phi(x,y) \leq e^{-\frac{|x|^2}{4} - \frac{|y-2H(x)\nu(x)|^2}{4}} \, \varphi(x)\] 
for all $(x,y) \in A$. Applying the change of variables formula to the map $\Phi$ gives 
\begin{align*} 
1 &= (4\pi)^{-\frac{n+1}{2}} \int_{\mathbb{R}^{n+1}} e^{-\frac{|\xi|^2}{4}} \, d\xi \\ 
&\leq (4\pi)^{-\frac{n+1}{2}} \int_\Sigma \bigg ( \int_{T_x^\perp \Sigma} e^{-\frac{|\Phi(x,y)|^2}{4}} \, |\det D\Phi(x,y)| \, 1_A(x,y) \, dy \bigg ) \, d\text{\rm vol}(x) \\ 
&\leq (4\pi)^{-\frac{n+1}{2}} \int_\Sigma \bigg ( \int_{T_x^\perp \Sigma} e^{-\frac{|x|^2}{4}-\frac{|y-2H(x)\nu(x)|^2}{4}} \, \varphi(x) \, dy \bigg ) \, d\text{\rm vol}(x) \\ 
&= (4\pi)^{-\frac{n}{2}} \int_\Sigma e^{-\frac{|x|^2}{4}} \, \varphi(x) \, d\text{\rm vol}(x).
\end{align*}
This shows that $\int_\Sigma \varphi \, d\gamma \geq 1$. To summarize, we know that 
\[\int_\Sigma \varphi \, \log \varphi \, d\gamma - \int_\Sigma \frac{|\nabla^\Sigma \varphi|^2}{\varphi} \, d\gamma - \int_\Sigma \varphi \, \Big ( H-\frac{1}{2} \, \langle x,\nu \rangle \Big )^2 \, d\gamma = 0\] 
and 
\[\bigg ( \int_\Sigma \varphi \, d\gamma \bigg ) \, \log \bigg ( \int_\Sigma \varphi \, d\gamma \bigg ) \geq 0.\] 
From this, the assertion follows. \\ 

Theorem \ref{log.Sobolev} is particularly useful on self-similar shrinking solutions to mean curvature flow. To explain this, suppose that $\Sigma$ is a hypersurface in $\mathbb{R}^{n+1}$. We say that $\Sigma$ is a self-similar shrinker if $H = \frac{1}{2} \, \langle x,\nu \rangle$. Self-similar shrinkers can be characterized as critical points of the Gaussian area $\gamma(\Sigma) = (4\pi)^{-\frac{n}{2}} \int_\Sigma e^{-\frac{|x|^2}{4}} \, d\text{\rm vol}$. The Gaussian area appears naturally in connection with Huisken's monotonicity formula for mean curvature flow \cite{Huisken}. Self-similar shrinkers achieve equality in Huisken's monotonicity formula; they play a central role in understanding singularity formation in mean curvature flow (see \cite{Colding-Minicozzi2}, \cite{Huisken}).

\begin{corollary}
\label{log.Sobolev.for.shrinkers}
Let $\Sigma$ be a compact hypersurface in $\mathbb{R}^{n+1}$ without boundary, and let 
\[d\gamma = (4\pi)^{-\frac{n}{2}} \, e^{-\frac{|x|^2}{4}} \, d\text{\rm vol}\] 
denote the Gaussian measure on $\Sigma$. If $\Sigma$ is a self-similar shrinker, then 
\[\int_\Sigma \varphi \, \log \varphi \, d\gamma - \int_\Sigma \frac{|\nabla^\Sigma \varphi|^2}{\varphi} \, d\gamma \leq \bigg ( \int_\Sigma \varphi \, d\gamma \bigg ) \, \log \bigg ( \int_\Sigma \varphi \, d\gamma \bigg )\]
for every positive smooth function $\varphi$ on $\Sigma$.
\end{corollary}

\end{document}